\documentclass[sn-basic]{sn-jnl}
\usepackage[utf8]{inputenc}
\usepackage{natbib}
\usepackage{geometry}
\RequirePackage{etex}
\usepackage{graphicx}%
\usepackage{multirow}%
\usepackage{amsmath,amssymb,amsfonts}%
\usepackage{amsthm}%
\usepackage{mathrsfs}%
\usepackage[title]{appendix}%
\usepackage{xcolor}%
\usepackage{textcomp}%
\usepackage{manyfoot}%
\usepackage{booktabs}%
\usepackage{algorithm}%
\usepackage{algorithmicx}
\usepackage[noend]{algpseudocode}
\usepackage{listings}%
\usepackage{comment}%
\usepackage{tikz}
\usepackage{tikz-network}
\usepackage{xfrac}
\usepackage{mathtools}
\usetikzlibrary{matrix}
\usepackage{bbm}
\usepackage[normalem]{ulem}
\usepackage{anyfontsize}
\usepackage{nicematrix}
\usepackage{blkarray}
\usepackage{bbm}
\usepackage[normalem]{ulem}

\usepackage[capitalize,nameinlink,compress]{cleveref} 
\crefname{figure}{Figure}{Figures} 
\crefname{equation}{}{} 
\crefname{assumption}{Assumption}{Assumptions}
\crefname{subsection}{Subsection}{Subsections}

\newtheorem{theorem}{Theorem}
\newtheorem{proposition}[theorem]{Proposition}%
\newtheorem{example}{Example}%
\newtheorem{remark}{Remark}%
\newtheorem{lemma}{Lemma}%
\newtheorem{corollary}{Corollary}%
\newtheorem{assumption}{Assumption}

\newtheorem{definition}{Definition}%

\newcommand{\R}{\mathbb{R}}

\newcommand{\TP}{\mathbb{TP}}

\newcommand{\one}{{\bf{1}}}
\newcommand{\zeros}{{\bf{0}}}

\DeclareMathOperator*{\argmin}{argmin}

\newcommand{\fw}{\mathrm{FW}}
\newcommand{\TPT}[1]{ \R^{#1}/\R \mathbf{1}}

\newcommand{\tr}{\mathrm{tr}}
\newcommand{\sample}{p_1, \dots, p_m}
\newcommand{\nset}{\mathcal{S}}

\newcommand{\CovDec}{\mathrm{CovDec}}
\newcommand{\BCD}{\mathrm{BCD}}
\newcommand{\mcf}{\mathrm{MCF}}
\newcommand{\ot}{\mathrm{OT}}
\newcommand{\program}[1]{(\mathrm{#1})}
\newcommand{\mydef}{\coloneq}
\newcommand{\D}{\mathrm{d}}
\newcommand{\dtr}{\mathrm{d}_{\mathrm{tr}}}
\newcommand{\drtr}{\mathrm{d}_{\gamma}}
\newcommand{\tnorm}[1]{\lvert\lvert #1 \rvert\rvert_{\tr}} 
\newcommand{\Rmax}{\mathbb{R}_{\max}}
\newcommand{\Rmin}{\mathbb{R}_{\min}}
\newcommand{\Xset}{\mathcal{X}}

\newcommand{\edit}[1]{\textcolor{black}{#1}}

\begin{document}
\title{Tropical Fermat--Weber Polytropes}


\author[1]{\fnm{John} \sur{Sabol}}\email{john.sabol@nps.edu}

\author[2]{\fnm{David} \sur{Barnhill}}\email{david.barnhill@nps.edu}

\author*[1]{\fnm{Ruriko} \sur{Yoshida}}\email{ryoshida@nps.edu}

\author[3]{\fnm{Keiji} \sur{Miura}}\email{miura@kwansei.ac.jp}

\affil*[1]{\orgdiv{Department of Operations Research}, \orgname{Naval Postgraduate School}, \orgaddress{\street{1 University Circle}, \city{Monterey}, \postcode{93943}, \state{CA}, \country{USA}}}

\affil[2]{\orgdiv{Department of Mathematics}, \orgname{U.S. Naval Academy}, \orgaddress{\street{121 Blake Road}, \city{Annapolis}, \postcode{21012}, \state{MD}, \country{USA}}}

\affil[3]{\orgdiv{Department of Biosciences}, \orgname{Kwansei Gakuin University}, \orgaddress{\street{1 Gakuen Uegahara}, \city{Sanda}, \postcode{669-1330}, \state{Hyogo}, \country{Japan}}}

\abstract{We study the geometry of tropical Fermat--Weber points, that is, optimal solutions to a location problem over a projective space using a dissimilarity measure derived from the tropical metric. It is well-known that for a given sample, such points are not necessarily unique, and we show that the set of all possible Fermat--Weber points forms a polytrope. This follows from the fact that our location problem turns out to be dual to a particular minimum-cost flow problem, and we describe the polytrope of optimal locations in the terminology of tropical geometry. We also provide a simple gradient descent algorithm that converges to the Fermat--Weber polytrope.}

\keywords{max-plus algebra, Fermat--Weber, min-cost flow}

\maketitle

\section{Introduction}\label{sec:Intro}

Let $\Xset$ be a set and $\D:\Xset \times \Xset \rightarrow \R$ be a bivariate function. For some finite sample $\nset \mydef \{\sample\}\subset \Xset$, we wish to find $x \in \Xset$ that minimizes $f(x)\mydef \sum_{i=1}^m\D(x,p_i)$. Such a minimizer, $x^*$, is called an \textit{M-estimator} in the statistical literature. If $\Xset$ is a compact topological space and $x\mapsto \D(x,p_i)$ is lower semi-continuous for each fixed $p_i$, then the existence of such a minimizer is guaranteed by standard results in optimization theory. When, in addition, $\D$ is a dissimilarity measure, i.e., $\D(x,y)\geq0$ and $\D(x,x)=0$, we call $x^*$ a \textit{generalized median}. A generalized median that is also symmetric (i.e., $\D$ is a metric), is called a \textit{Fr\'echet median}. Finally, when $\Xset$ is a normed real vector space, we call $x^*$ a \emph{geometric median}.

In some applications it is desirable to associate a weight, $w_i\in \R_{\geq0}$, to each sample element $p_i\in \nset$. As long as these weights are non-negative, existence is not affected. Thus, our minimizers may take the more general form
\[
x^* = \argmin_{x\in \Xset} f(x) \qquad \text{where} \qquad f(x) = \sum_{i=1}^m w_i\D(x, p_i).
\]
The problem of finding a (weighted) geometric median is often called the {\em Fermat--Weber location problem}, and any optimal point $x^*$ a {\em Fermat--Weber point} of the sample. In what follows we use $\Xset^*$ to denote the set of all such minimizers, i.e. $\Xset^* \mydef \{x \in \Xset \mid f(x) \leq f(y),\, y\in \Xset \}$.

Since \cite{Weiszfeld}, substantial work has been done on geometric medians and their weighted extensions (e.g., see \citep{Brimberg} and references therein). In contrast, considerably less attention has been paid to non-Euclidean settings, where many important classes of problems naturally arise. Properly addressing these problems requires statistical techniques specifically adapted to the geometry of the underlying space. For example, in phylogenetics, $x^*$ is a Bayes estimator if $\nset$ is drawn from the posterior distribution over the space of phylogenetic trees \citep{bayesEstimator}. In other areas, observations may lie on curved manifolds, such as the space of positive definite matrices, where standard Euclidean methods fail to capture the intrinsic geometry. Identifying efficient means for computing statistical measures over such spaces is an important endeavor.

The space we consider in this paper is the tropical projective torus, $\Xset=\TPT{n}$ ($n\geq 2$), a quotient space that provides a natural setting for tropical geometry. We will properly define $\TPT{n}$ later, deferring all definitions to \cref{sec:basics}. 

Research on tropical analogues of geometric medians was initiated in \citep{LSTY, Lin2016TropicalFP}, where the authors applied the tropical metric $\dtr$ to uniformly-weighted samples from $\TPT{n}$. They showed that the set $\Xset^*$ of \emph{tropical Fermat--Weber points} forms a classically convex polytope in $\R^{n-1}$ (which is isomorphic to $\TPT{n}$), and provided a linear program to compute such a minimizer. In \citep{Com_neci_2023, Comneci2024}, the authors introduce a quasi-metric counterpart to $\dtr$ which they call the \emph{asymmetric tropical distance}, denoted by $\D_\Delta$. Applying $\D_\Delta$ to samples from $\TPT{n}$, they show that an \emph{asymmetric Fermat--Weber point} $x^*$ corresponds to an optimal solution of a particular transportation problem. In \citep{cox2023}, the authors extend such results by showing that (positively) weighted asymmetric tropical Fermat--Weber points retain the same tropical convexity properties as their unweighted counterparts.

In this paper, we extend the previous work to a more general setting by considering a quasi-metric on $\TPT{n}$ that, in essence, provides a weighted interpolation between $\dtr$ and $\D_{\Delta}$. We revisit the tropical Fermat--Weber problem over $\TPT{n}$ using this measure. A summary of our contributions are as follows:
\begin{itemize}
    \item We show that this more general formulation of the weighted tropical Fermat--Weber problem is solved by a minimum-cost flow problem which has a corresponding optimal transport formulation.
    \item We prove that $\Xset^*$ forms a bounded cell in the common refinement given by max-plus and min-plus covector decompositions. This allows us to completely characterize $\Xset^*$ in terms of its extreme vertices and implies efficient methods for recovering \emph{all} tropical Fermat--Weber solutions for a given problem.
    \item We formulate a tropical Fermat--Weber (projective) subgradient descent algorithm and provide the geometric intuition underpinning several practical improvements.
\end{itemize}

The organization of this paper is as follows. In \cref{sec:basics} we provide a brief overview of the background necessary to connect several tropical geometry and network optimization concepts. We utilize \citep{joswigBook} and \citep{AMO_NetworkFlows} as our primary sources for tropical geometry and network optimization respectively. In \cref{sec:Compute_FW}, we provide the min-cost flow and optimal transport formulations for computing $\Xset^*$. In \cref{sec:subgradient}, we consider the subdifferential of the Fermat--Weber objective function and discuss the geometric interpretation behind an iterative subgradient descent algorithm over $\TPT{n}$. 
We conclude in \cref{sec:conclusion} with a discussion of future research directions. 

Regarding notation, we use $[m]$ to denote the set $\{1,\ldots,m\}$ for positive integer $m$. All vectors are considered column vectors, and $(\cdot)^{\top}$ denotes the transpose of a vector/matrix. We let $\one\mydef(1,\ldots,1)^{\top}$ and $\boldsymbol{0}\mydef(0,\ldots,0)^{\top}$ represent the vector of all ones and all zeros respectively. $\R\one \mydef \{c \one \in \R^n\mid c \in \R\}$ can then be seen to represent any vector with common elements. Finally, $\emptyset$ represents the empty set.

\section{Preliminaries}\label{sec:basics}
\textbf{Tropical Geometry Part I.}
Tropical geometry is the geometric counterpart to tropical algebra, an algebra that replaces addition with minimization and multiplication with addition. That is, for $a,\,b\in\R$, $a\oplus b \mydef \min(a,b)$ and  $a \odot b \mydef a+b$ where $\oplus$ and $\odot$ are used to distinguish these ``tropical'' operators from their classical analogues. Under the tropical algebra, $0$ is the multiplicative identity element ($a\odot 0=a$). To obtain an additive identity we extend $\R$ by defining $\Rmin \mydef \R\cup \{\infty\}$, thereby allowing $\infty$ to assume this role ($a\oplus\infty=a$). Together, $(\Rmin, \oplus, \odot)$ defines the ``min-plus'' semiring, also called the \emph{tropical semiring}. 

Operations in the tropical semiring extend naturally to scalar multiplication and vector addition through component-wise application, yielding the \emph{tropical semimodule} $\Rmin^n$. Matrix multiplication proceeds as expected. For example, we have
\[
A=\begin{pmatrix}
    0 & 4 \\
    \infty & -2 \\
    5 & 1
\end{pmatrix},\qquad
B=\begin{pmatrix}
    1 & 3 & \infty\\
    2 & 4 & 0
\end{pmatrix}, \quad
A\odot B=\begin{pmatrix}
    1 & 3 & 4 \\
    0 & 2 & -2 \\
    3 & 5 & 1
\end{pmatrix}.
\]
Since this paper also uses standard notation from linear programming (e.g., $Ax=b)$, we annotate tropical operators explicitly throughout.

\emph{Tropical projective space} is the quotient space, $\TP^{n-1}\mydef(\Rmin^n\setminus\{\infty\one\})/\R\one$, obtained by identifying vectors that differ by tropical scalar multiplication. As in classical projective geometry, we remove the zero vector (here, $\infty\one$) prior to quotienting. Since tropical scalar multiplication corresponds to adding a constant to every coordinate, two vectors $x,\,y\in \TP^{n-1}$ represent the same projective point precisely when $x=y+\lambda\one$ for some $\lambda\in \R$.

The \emph{tropical projective torus}, $\TPT{n}$, is the subset of points in $\TP^{n-1}$ with finite coordinates. Since each element $x \in \TPT{n}$ is an equivalence class under the relation $x \sim x+\lambda\one$, it is common to select a representative from each equivalence class by restricting to an $(n-1)$-dimensional hyperplane, e.g. $\{x\in \R^n\mid x_j=0\}$ for some $j\in[n]$ or $\{x\in \R^n\mid \sum_jx_j=0\}$. Under this identification, $\TPT{n}$ is homeomorphic to $\R^{n-1}$ and therefore carries the standard Euclidean topology. This finite-dimensional linear structure motivates our use of the term ``Fermat--Weber'' problem rather than the more general terms Fr\'echet median or generalized median.

$\TPT{n}$ is a metric space equipped with a norm $\tnorm{x}\mydef \max_jx_j - \min_jx_j$ and its associated metric $\dtr(\cdot,\cdot)$, which we call the \emph{tropical distance}. This distance appears in several areas of mathematics under various names, including the \emph{Hilbert projective metric} and the \emph{span semi-norm}. Strictly speaking, the Hilbert projective metric $d_H$ is distinct from $\dtr$. The former arises as a pseudometric over a closed convex cone $K$ in a Banach space $\mathcal{B}$, where $K$ is pointed, i.e. $K\cap (-K)=\{\boldsymbol{0}\}$. Thus, $K$ determines a partial order $\leq_K$ on $\mathcal{B}$ and we can define $M(u/v)\mydef \inf\{\lambda\mid u \leq_K \lambda v\}$ and
$m(u/v)\mydef \sup\{\mu \mid \mu v \leq_K u \}$
for $u,v \in K\setminus\{\boldsymbol{0}\}$. $d_H$ is then defined as
\[
\D_H(u,v)\mydef \log\left( \frac{M(u/v)}{m(u/v)}\right).
\]
A key feature of $\D_H$ is its invariance under multiplication of either argument by a positive scalar. Consequently, it depends only on the rays generated by $u$ and $v$, and therefore induces a metric on the set of rays of $K$, i.e., on the projectivization of the cone. In the special case where $K=\R^n_{>0}$ is the positive orthant, the extrema defining $M$ and $m$ are attained, ensuring that $\D_H(u,v)$ is finite. The resulting metric on the projective space associated with $K$ is the standard form of Hilbert's projective metric. By taking $x,\,y$ such that $x_j\mydef \log(u_j)$ and $y_j\mydef \log(v_j)$, we have through simple substitution, that $\D_H(u,v)=\dtr(x,y)$. Thus, in a certain sense we are justified in thinking of $\TPT{n}$ as $\R^n_{>0}$ viewed through ``logarithmic'' glasses.

An asymmetric variant of tropical distance was introduced in \citep{Com_neci_2023}, which the authors call, the \emph{asymmetric tropical distance} $\D_{\Delta}(\cdot,\cdot)$. In summary, we have
\begin{align*}
    \dtr(x,y) &\mydef \max_{j\in[n]}(x_j-y_j)-\min_{j\in[n]}(x_j-y_j),\\
    \D_\Delta(x,y) &\mydef \max_{j\in[n]}(x_j - y_j) - \frac{1}{n}\sum_{j=1}^n(x_j - y_j).
\end{align*}
The next definition essentially merges these two distances. 

\begin{definition}[Generalized Tropical Distance]
    For $x,y\in \TPT{n}$ and $\gamma \in [0,1]$, the \emph{generalized tropical distance} between $x$ and $y$ is
    \begin{equation}\label{eq:trop_dist}
    \drtr(x,y) \mydef \gamma\max_{j\in[n]}(x_j-y_j) - (1-\gamma)\min_{j\in[n]}(x_j-y_j) + \frac{(1-2\gamma)}{n}\sum_{j=1}^n(x_j - y_j).
\end{equation}
\end{definition}

Intuitively, $\gamma$ can be thought of as a direction ratio governing the relative importance of the one-way distances $x\rightarrow y$ vs. $y\rightarrow x$. It is easy to see that $\D_{1}(x,y) = \D_\Delta(x,y)$, $\D_{0}(x,y) = \D_\Delta(y,x)$, and $\D_{\gamma}(x,y)+\D_{1-\gamma}(x,y)=\dtr(x,y)$. Like its specialized counterparts, $\drtr$ is non-negative and separates equivalence classes of $\TPT{n}$, i.e. it is invariant under scalar shifts: $\drtr(x,y)=\drtr(x+a\one,y+b\one)$ for $a,b\in \R$, and $\drtr(x,y)=0$ \emph{if and only if} $x=y+b\one$ for some $b\in \R$. In general, $\drtr$ is asymmetric, $\gamma=1/2$ being the exception. For fixed $y$, $\drtr(x,y)$ is a continuous and piecewise linear function. {The next lemma shows it is also coercive, which, since $\Xset=\TPT{n}$ is finite-dimensional, means it attains its minimum.  
\begin{lemma}\label{lem:coercive}
    A generalized tropical distance is coercive in $\TPT{n}$.
\end{lemma}
\begin{proof}
    Since $\drtr(x,y)=\drtr(x-y,\zeros)$, we may assume $y=\zeros$ without loss of generality. It suffices to show that $\drtr(x,\zeros)\geq \tfrac{1}{n}\lVert x \rVert_{\mathrm{tr}}$. Let $x^{(i)}$ denote the $i$-th order statistic of $x$ after taking an arbitrary fixed coordinate representative. Observe that,
    \begin{align*}
    \D_0(x,\zeros)&= \tfrac{1}{n}\sum_{i=1}^n\bigl(x^{(i)}-x^{(1)}\bigr)=\tfrac{1}{n}\bigl(x^{(n)}-x^{(1)}\bigr)+\tfrac{1}{n}\sum_{i=1}^{n-1}\bigl(x^{(i)}-x^{(1)}\bigr)\geq \tfrac{1}{n}\lVert x\rVert_{\mathrm{tr}},\\
    \D_1(x,\zeros)&= \tfrac{1}{n}\sum_{i=1}^n\bigl(x^{(n)}-x^{(i)}\bigr)=\tfrac{1}{n}\bigl(x^{(n)}-x^{(1)}\bigr)+\tfrac{1}{n}\sum_{i=2}^{n}\bigl(x^{(n)}-x^{(i)}\bigr)\geq \tfrac{1}{n}\lVert x\rVert_{\mathrm{tr}}.
    \end{align*}
    Thus, $\drtr(x,\zeros)=\gamma \mathrm{d}_0(x,\zeros)+(1-\gamma)\mathrm{d}_1(x,\zeros)\geq \tfrac{\gamma}{n}\lVert x\rVert_{\mathrm{tr}} + \tfrac{1-\gamma}{n}\lVert x\rVert_{\mathrm{tr}}=\tfrac{1}{n}\lVert x\rVert_{\mathrm{tr}}$.
\end{proof}}

We are now ready to define our tropical Fermat--Weber objective function. Consider a sample $\nset=\{\sample\}\subset\TPT{n}$, with each point given a direction ratio $\gamma_i\in[0,1]$ and weight $w_i\in\R_{\geq0}$. Let $\gamma=(\gamma_1,\ldots,\gamma_m)$ and $w=(w_1,\ldots,w_m)$. We assume that $m \geq n$ (it is conceivable that in many applications $m \gg n$). The tropical Fermat--Weber objective function, $f_\nset(x)$, is then defined as
\begin{equation}\label{eq:trop_fw_obj}
    f_\nset(x) \mydef \sum_{i=1}^mw_i\D_{\gamma_i}(x,p_i),
\end{equation}
where $\D_\gamma(x,\cdot)$ is given by \cref{eq:trop_dist}. \edit{Since positively weighted sums of coercive functions is again coercive, we see that the optimal solution set $\Xset^*\neq \emptyset$ to \cref{eq:trop_fw_obj} is non-empty.}

Recall that a set $C\subset\R^n$ is (classically) convex if $\lambda x + \mu y \in C$ for all $x,y\in C$ and $\lambda,\mu\in \R_{\geq0}$ such that $\lambda+\mu=1$. A tropical analogue to convexity was introduced in \citep{DS} under the name \emph{tropical convexity}.
\begin{definition}[Tropically Convex Set]\label{def:tropical_convex_set}
    A subset $C\subset \TP^{n-1}$ is \emph{tropically convex} if, for all $x,y\in C$ and $\lambda,\mu \in \Rmin$ such that $\lambda\oplus \mu =0$, we have $\lambda \odot x \oplus \mu \odot y \in C$.
\end{definition}
When $C\subset \TPT{n}$ (and $\lambda,\mu \in \R$) we can drop the requirement that $\lambda\oplus \mu =0$, which then matches the definition in \citep{DS}.
For a subset $C \subset \Rmin^n$, the \emph{tropical convex hull}
\[
\mathrm{tconv}(C) \mydef \{\lambda_1 \odot c_1 \oplus \dots \oplus \lambda_m \odot c_m \mid \lambda_i \in \R, c_i \in C\}
\]
is the smallest tropically convex subset containing $C$. The set $\mathrm{tconv}(C)$ is called a \emph{tropical polytope} if it admits a finite generating set. A tropical version of Carath\'eodory's Theorem from classical convexity states that any tropical polytope $\mathrm{tconv}(C) \subset \Rmin^n$ can be generated by at most $n$ extreme points of $\mathrm{tconv}(C)$. We call these extreme points the \emph{tropical vertices} of the polytope to distinguish them from (ordinary) vertices. For a proof of the theorem, see \citep{DS,joswigBook}.

\begin{remark}
    Classic convexity and tropical convexity are, in general, distinct. A tropical polytope that is also classically convex is called a \emph{polytrope}.
\end{remark}

Going back to the beginning of \Cref{sec:basics} and replacing ``$\min$'' with ``$\max$'' in the definitions yields the \emph{max-plus} semiring, $(\Rmax, \boxplus, \odot)$, where $\Rmax \mydef \R\cup\{-\infty\}$ and $a\boxplus b\mydef \max(a,b)$. The isomorphism between min-plus and max-plus semirings is reflective around $x=\boldsymbol{0}$, as demonstrated by the identity $\max(a,b)=-\min(-a,-b)$. Min-plus convexity does not imply max-plus convexity, nor the reverse. A set that is both min-plus convex and max-plus convex is called \emph{bi-tropically convex}. It is easy to show bi-tropical convexity implies classical convexity and hence, any \edit{closed} bi-tropical set is a polytrope, see \citep[Chapter 5]{joswigBook}.

In general, the set of tropical vertices of a tropical polytope can be arbitrarily large. However, this situation improves dramatically for polytropes. In fact, for an $n$-dimensional polytrope $\mathcal{P}$:

\begin{enumerate}
    \item $\mathcal{P}$ is a tropical simplex with $n+1$ tropical vertices \citep[Lemma 10]{DS}.
    \item $\mathcal{P}$ has at most $n(n+1)$ facets, and this bound is sharp \citep[Proposition 5]{JoswigKulas_Trop_Ordinary_Combined}.
    \item $\mathcal{P}$ has at most $\binom{2n}{n}$ vertices, and this bound is sharp \citep[Proposition 19]{DS}; \citep[Theorem 2.3(2)]{Gelfand_1997_Combinatorics}.
\end{enumerate}

Recall that for a convex set $X$, a function $f:X \rightarrow\R$ is called convex if its epigraph, $\mathrm{epi}(f)\mydef \{(x,t) \in X \times \R\mid f(x)\leq t\}$, is convex. Equivalently, $f$ is convex if and only if $f(\lambda x + \mu y) \leq \lambda f(x) + \mu f(y)$ for all $x,y\in X$, $\lambda,\mu \in \R_{\geq0}$ such that $\lambda+\mu=1$. 

\begin{proposition}\label{prop:f_is_convex}
    $f_\nset(x)$ as defined in \Cref{eq:trop_fw_obj} is convex.
\end{proposition}
\begin{proof}
    It suffices to show that $w_i\D_{\gamma_i}(x,p_i)$ is convex for each $i\in [m]$. Index $i$ is therefore arbitrary and we suppress it in what follows. Let $\alpha\mydef w\gamma$ and $\beta \mydef w(1-\gamma)$.
    \begin{align*}
    w\D_{\gamma}(x,p)&=\alpha\max_j(x_j-p_j)-\beta\min_j(x_j-p_j)+(\beta-\alpha)\tfrac{1}{n}\sum_j(x_j-p_j),\\
    &=\alpha \, \Phi_p(x)-\beta \, \Psi_p(x)+(\beta-\alpha)\tfrac{1}{n}\sum_j(x_j-p_j),\\
    &=\alpha \, \Phi_p(x)+\beta \overline{\Psi}_p(x) + (\beta-\alpha)\tfrac{1}{n}\sum_jx_j + (\alpha-\beta)\tfrac{1}{n}\sum_jp_j,
    \end{align*}
    where $\Phi_p(x) \mydef \max_j(x_j - p_j)$ and $\overline{\Psi}_p(x) = -\Psi_p(x) \mydef \max_j(p_j-x_j)$. Both $\Phi_p(x)$ and $\overline{\Psi}_p(x)$ are maximums over affine functions and are therefore convex in $x$. Multiplication by non-negative scalars preserves convexity. Thus, $\D_{\gamma}(x,p)$ is a sum of convex functions. A (positively weighted) sum of convex functions is convex.
\end{proof}

The tropical analogue of convex function is a \emph{tropically convex function}.
\begin{definition}[Tropically Convex Function]\label{def:maxplus_convex_function}
    Let $\overline{\R}\mydef \R\pm\{\infty \}$. A map $f:\Rmin^n \rightarrow \overline{\R}$ is tropically convex if its epigraph is tropically convex. Equivalently, $f$ is tropically convex if, and only if, for all $x,y\in \Rmin^n$, and $\lambda,\mu \in \Rmin$ such that $\lambda \oplus \mu= 0$, $f(\lambda \odot x \oplus \mu\odot y) \leq \lambda \odot f(x) \oplus \mu \odot f(y)$.
\end{definition}

\begin{proposition}
    $f_\nset(x)$ as defined in \Cref{eq:trop_fw_obj} is generally neither min-plus nor max-plus convex.
\end{proposition}
\begin{proof}
    We provide the following simple counter-example for $m=3$, $n=3$.
    \begin{figure}
        \centering
        \begin{minipage}{0.7\textwidth}
            \centering
            {\small
            \begin{align*}
            x_1 &= (0,0,-1)^{\top}, \quad y_1 = (0,-1,0)^{\top}, \quad & f_\nset(x_1) &= f_\nset(y_1)=2,\\
            x_2 &= (0,\tfrac{3}{2},1)^{\top}, \quad ~~y_2 = (0,1,\tfrac{3}{2})^{\top}, \quad & f_\nset(x_2) &= f_\nset(y_2)=2,\\
            z_1 &= \lambda\odot x_1 \oplus \mu\odot y_1 = (0,-1,-1)^{\top}, \quad & f_\nset(z_1) &= \tfrac{9}{4},\\
            z_2 &= \lambda\odot x_2 \boxplus \mu\odot y_2 = (0,\tfrac{3}{2},\tfrac{3}{2})^{\top}, \quad & f_\nset(z_2) &= \tfrac{5}{2}.
            \end{align*}
            }
        \end{minipage}
        \hfill
        \begin{minipage}{0.28\textwidth}
            \centering
            \begin{tikzpicture}
            \draw[step=1cm, gray!30, thin] (-1,-1) grid (2,2);
            \draw[-,gray!60] (-1,0) -- (2,0);
            \draw[-,gray!60] (0,-1) -- (0,2);
            
            \draw (-1,0) -- (-1/3,1) -- (0.05,1.35) -- (0.5,1.5) -- (1,1.5) -- (1.06,1.47) -- (1.3,1.35) -- (1.35,1.3) -- (1.37,1.26) -- (1.5,1) -- (1.5,0.5) -- (1.35,0.05) -- (1,-1/3) -- (0,-1) -- (-2/3,-2/3) -- cycle;
            \draw [-,dotted](0,-1) -- (-1,-1) -- (-1,0);
            \draw [-,dotted](3/2,1) -- (3/2,3/2) -- (1,3/2);
            
            \node at (0,0)[circle,fill,inner sep=1.5pt]{};
            \node at (0,1)[circle,fill,inner sep=1.5pt]{};
            \node at (1,0)[circle,fill,inner sep=1.5pt]{};
            \node at (-1,0)[circle,draw=black,fill=white,inner sep=1.5pt]{};
            \node at (0,-1)[circle,draw=black,fill=white,inner sep=1.5pt]{};
            \node at (-1,-1)[circle,draw=black,fill=white,inner sep=1.5pt]{};
            \node at (3/2,1)[circle,draw=black,fill=white,inner sep=1.5pt]{};
            \node at (1,3/2)[circle,draw=black,fill=white,inner sep=1.5pt]{};
            \node at (3/2,3/2)[circle,draw=black,fill=white,inner sep=1.5pt]{};
            \node[] at (-0.2,0.2) {\small$p_1$};
            \node[] at (0.25,0.8) {\small$p_3$};
            \node[] at (0.8,0.2) {\small$p_2$};
            \node[] at (0,-1.3) {\small$x_1$};
            \node[] at (-1.25,0) {\small$y_1$};
            \node[] at (-1.25,-1) {\small$z_1$};
            
            \node[] at (1.8,0.8) {\small$x_2$};
            \node[] at (0.8,1.7) {\small$y_2$};
            \node[] at (1.7,1.7) {\small$z_2$};
            \end{tikzpicture}
        \end{minipage}
        \caption{Let $\nset=\{p_1,p_2,p_3 \}$ with $p_1=(0,0,0)^{\top}$, $p_2=(0,1,0)^{\top}$, and \edit{$p_3=(0,0,1)^{\top}$}. We take the unweighted, symmetric case given by $w=(1,1,1)$ and $\gamma=(\tfrac{1}{2},\tfrac{1}{2},\tfrac{1}{2})$, and take $\lambda=\mu=0$. The black solid line depicts the $f_\nset(\cdot)=2$ level set. 
        We see that $f_\nset(z_1)=9/4$ and $f_\nset(z_2)=5/2$. Thus, we conclude that $f_\nset(\cdot)$ is neither min-plus convex nor max-plus convex at this level set, and \edit{therefore} is not a min-plus or max-plus convex function.}
        \label{fig:counterexample_trop_convexity}
    \end{figure}
\end{proof}

\textbf{Network Optimization.} We transition to briefly review concepts from network optimization and graph theory, returning to tropical geometry to highlight explicit connections as they arise. The following definitions can be found in any standard text.

Let $\mathcal{G}(V,E)$ denote a directed graph 
where $V$ is a finite set of nodes and $E$ is a finite set of arcs. Each arc $e\in E$ is an ordered pair $e=(i,j)\in V\times V$ where $i$ is the tail and $j$ is the head. A self-loop is an arc of the form $(i,i)$. Two distinct arcs, $e_1=(i,j)$ and $e_2=(i,j)$, are called parallel if they share the same tail and head. $\mathcal{G}$ is called simple if it contains no self-loops and no parallel arcs. For any $i_0\in V$, a \emph{walk} along $\mathcal{G}$ is a sequence of nodes $\{i_0,i_1,\ldots,i_L\}\subseteq V$ such that every adjacent pair in the sequence corresponds to an arc in $\mathcal{G}$, i.e. $(i_{\ell-1},i_\ell)\in E$ for every $\ell \in [L]$. Equivalently, a walk can be described as an ordered sequence of arcs \edit{$(i_{0},i_1),\ldots,(i_{\ell-1},i_\ell)$}. A walk is \emph{simple} if no node is repeated. Simple walks are called \emph{paths}. Node $j$ is \emph{reachable} from $i$ if there exists a path from $i$ to $j$. $\mathcal{G}$ is \emph{strongly connected} if every node is reachable from every other node. $\mathcal{G}$ is \emph{weakly connected} if the underlying undirected graph obtained by ignoring arc directions is connected. A \emph{cycle} is a closed walk in which no node is repeated except the first and last. The \emph{weight} of a walk/path/cycle is the sum of arc costs over the arcs in the sequence (counted with multiplicity).

Take $\mathcal{G}(V,E)$, $E\subset V \times V$, to be a simple, directed graph with $|V|=N$ nodes and $|E|=M$ arcs. This assumption is without loss of generality, as any directed graph with loops or parallel arcs can be transformed into a simple graph by introducing auxiliary nodes and splitting arcs, see e.g. \citep{AMO_NetworkFlows}. Let $b\in\R^N$ denote the vector of node demands, and let $c\in\R^M$ denote the vector of arcs costs. We call an $i\in V$ a \emph{supply node} if $b_i<0$, a \emph{demand node} if $b_i>0$, and a \emph{transhipment node} if $b_i=0$. The uncapacitated \emph{min-cost flow} ($\mcf$) problem with flow variables $x=\{x_{ij}\mid (i,j)\in E\}$ can be formulated as the linear program
\[
\mathrm{(MCF)}\qquad \qquad \min_{x
\geq0} c^{\top}x \quad \text{subject to} \quad Ax=b,
\]
where $A$ is the $N\times M$ matrix that has, for each arc $(i,j) \in E$, a column equal to $e_i - e_j$, with $e_k$ denoting the $k$-th standard basis vector of $\mathbb{R}^N$. We make the following standard assumptions on $\mcf$.
\begin{assumption}\label{assum:G_feasible}
    $\mathcal{G}$ is weakly connected, $\sum_ib_i=0$, and every demand node is reachable from at least one supply node.
\end{assumption}
\begin{assumption}\label{assum:G_bounded}
    $\mathcal{G}$ contains no cycle with negative weight. 
\end{assumption}
\cref{assum:G_feasible} is necessary for $\mcf$ to be feasible. \cref{assum:G_bounded} is necessary for $\mcf$ to have a bounded solution. The dual associated to $\mcf$ has a feasible region, $\mathcal{D}\subset \R^N$; a polyhedron defined by $A^{\top}\pi\leq c$ where each inequality is of the form $\pi_i-\pi_j \leq c_{ij}$. Polyhedra of this form are called \emph{weighted digraph polyhedra}, see \citep{joswigBook}. Thus, $\mathcal{D}$ is a weighted digraph polyhedron. {The next lemma, which is proved in \cite[Proposition 5.30]{joswigBook}, shows that any weighted digraph polyhedron is bi-tropically convex. We include a slightly modified version here for completeness.} 
\begin{lemma}[\cite{joswigBook}, Proposition 5.30]
    A weighted digraph polyhedron $\mathcal{D}$ is bi-tropically convex.
\end{lemma}\label{lem:WDP_are_tropConvex}
\begin{proof}
    We begin by proving the following inequality, which holds for any $a,b,c,d \in \R$.
    \[
    f\mydef \min(a,b)-\min(c,d)\leq \max(a-c,b-d)
    \]
    As there are only four possible cases, we enumerate them.
    \begin{enumerate}[leftmargin=2em]
        \item $a\leq b, c \leq d$: \quad $f=a-c\leq \max(a-c,b-d)$.
        \item $a\leq b, c > d$: \quad$f=a-d \leq \edit{b-d} \leq \max(a-c,b-d)$.
        \item $a > b, c\leq d$: \quad$f=b-c < \edit{a-c} \leq \max(a-c,b-d)$.
        \item $a > b, c > d$: \quad$f=b-d \leq \max(a-c,b-d)$.
    \end{enumerate}
    Thus, the inequality is proven. To prove the lemma, we need only show that $\mathcal{D}$ is closed under tropical addition and scalar multiplication. We assume $\mathcal{D} \neq \emptyset$, otherwise there is nothing to prove. Take arbitrary $\pi,\phi \in \mathcal{D}$, so that (by definition) $\pi_i-\pi_j\leq c_{ij}$ and $\phi_i-\phi_j\leq c_{ij}$ for every constraint in $\mathcal{D}$. Define $\psi \mydef \pi \oplus \phi$ and consider arbitrary $(i,j)$. We have
    \[
    \psi_i - \psi_j =\min(\pi_i,\phi_i)-\min(\pi_j,\phi_j) \leq \max(\pi_i - \pi_j, \phi_i-\phi_j) \leq c_{ij}.
    \]
    This shows that $\mathcal{D}$ is closed under tropical addition. Now take $\pi'\mydef \pi+c\one$ for $c\in\R$ and see that $\pi'_i-\pi'_j\leq c_{ij}$ for all $(i,j)$. This shows that $\mathcal{D}$ is closed under tropical scalar multiplication. Thus, $\mathcal{D}$ is min-plus convex. For max-plus, simply take
    \[
    \max(\pi_i,\phi_i)-\max(\pi_j,\phi_j)=\min(-\pi_j,-\phi_j)-\min(-\pi_i,-\phi_i)
    \]
    and apply the inequality as before.
\end{proof}

It is often convenient to represent a weighted digraph polyhedron $\mathcal{D}$ as a square matrix $C=(c_{ij})\in\Rmin^{n\times n}$ where each entry encodes the (non-redundant) inequality $\pi_i-\pi_j\leq c_{ij}$. To ensure every element of $C$ is well-defined, we can implicitly add the trivial inequalities $\pi_i-\pi_i\leq 0$ and $\pi_i-\pi_j \leq \infty$, which doesn't affect the geometry of $\mathcal{D}$.

In practical terms, $C$ constructed this way is equivalent to the weighted adjacency matrix of a directed graph on $n$ nodes. It is well-known that when this graph has no negative cycles, the path closure over the graph, $C^*$, is well-defined and computable using an all-pairs shortest path algorithm, e.g. Floyd-Warshall or Johnson's algorithm. The fact that this is a closure operation is best seen tropically via
\begin{equation}\label{eq:KleeneStar}
    C^*=\mathbb{I}_n\oplus C \oplus C^{\odot2} \oplus \cdots=\mathbb{I}_n\oplus C \oplus C^{\odot2} \oplus \cdots\oplus C^{\odot(n-1)}
\end{equation}
where $\mathbb{I}_n$ is the $n\times n$ tropical identity matrix (the square matrix with all entries equal to $\infty$ except for zeros along the diagonal) and $A^{\odot k}$ is the matrix $A$ ``powered up'' $k$ times (i.e., $A\odot \dots \odot A$). We can interpret $\mathbb{I}_n \oplus C$ as the step that incorporates the trivial inequalities described earlier. We call $C^*$ the \emph{Kleene star} of $C$ and note that the formula in \cref{eq:KleeneStar} is equivalent to the Bellman-Ford algorithm.

\cite{Sergeev2007_closures} examined such closures, noting that the columns of $C^*$ are precisely the (max-plus) vertices of the associated columns space, and \cite{JoswigKulas_Trop_Ordinary_Combined} showed how to compute these extreme (tropical) vertices from an ordinary inequality description (i.e., a weighted digraph polyhedron) using any all-pairs shortest path algorithm on a directed graph.

Assume the problem $\mcf$ is feasible and consider any optimal flow $x^*$ (which by our earlier assumptions must be finite). Let $\mathcal{D}^*\subseteq \mathcal{D}$ denote the space of optimal dual solutions. Strong duality guarantees $\mathcal{D}^*\neq \emptyset$ and so $\mathcal{D}^*$ is precisely the face of $\mathcal{D}$ obtained taking every positive optimal flow $\{x^*_{ij}\mid x^*_{ij}>0\}$ and changing the corresponding inequality in $\pi_i-\pi_j \leq c_{ij}$ into an equality. The next lemma shows that $\mathcal{D}^*$ is again a weighted digraph polyhedron.

\begin{lemma}[\cite{joswigBook}, Lemma 3.44]\label{lem:faces_of_WDP}
    Faces of weighted digraph polyhedra are weighted digraph polyhedra.
\end{lemma}
\begin{proof}
    Let $\mathcal{D}$ be an arbitrary weighted digraph polyhedron whose constraint set is indexed according to the arc set $E$ of its associated digraph, i.e. $\mathcal{D}\mydef \{\pi_i - \pi_j\leq c_{ij}\mid (i,j)\in E \}$. Let $E'\subseteq E$ be an arbitrary subset of arcs. For every arc $(i,j)\in E'$, take the corresponding inequality in $\mathcal{D}$, and turn it into an equality by adding the constraint $\pi_i-\pi_j \geq c_{ij}$. This larger constraint set defines a face of $\mathcal{D}$, which we denote by $\mathcal{D}'$. Rearranging the added constraints puts $\mathcal{D}'$ into the form of a weighted digraph polyhedron, which completes the proof.
\end{proof}

Recall that for a given a feasible flow $x$, the residual graph, denoted by $\mathcal{G}_x$, is the graph constructed from $\mathcal{G}$ by adding backwards arcs $(j,i)$ with negated costs $-c_{ij}$ for each arc $(i,j)$ that carries a positive flow.
Intuitively, \cref{lem:faces_of_WDP} describes the construction of the weighted digraph polyhedron corresponding to the residual graph of a feasible flow.
\edit{If $E'$ denotes a subset of arcs whose corresponding constraints in $\mathcal{D}$ are set to equality, then the subgraph $\mathcal{G'}=(V,E')$ is called an \emph{equality subgraph} of $\mathcal{G}$. If $x^*$ is an optimal flow in a network flow problem (e.g., $\mcf$), then complementary slackness enforces equality on every constraint $\pi_i - \pi_j\leq c_{ij}$ where $x^*_{ij}>0$. Hence, $x^*$ induces an equality subgraph, which we denote by $\mathcal{G}|_{x^*}=(V,E^*)$. Since $E^*$ encodes the ``tight'' constraints, we use $A^{T}_{E^*}$ to denote the sub-matrix of $A^{\top}$ whose rows satisfy $A^{\top}_{E^*}\pi=c_{E^*}$ in $\mcf$.} 

In general, the weighted digraph polyhedron corresponding to a residual graph may be empty due to the addition of incompatible constraints. This is equivalent to the introduction of a negative cycle in the residual graph. Fortunately, the complementary slackness conditions at optimality preclude the existence of such negative cycles in $\mathcal{G}_{x^*}$, a classical result in network optimization theory.

This means that we can take the path closure of the matrix associated with $\mathcal{D}^*$ (by computing the all-pairs shortest paths on $\mathcal{G}_{x^*}$) and thereby obtain all tropical extreme points of $\mathcal{D}^*$. 

\begin{lemma}\label{lem:MCF_duals_are_tropConvex}
    The set of feasible dual solutions $\mathcal{D}$ and the set of optimal dual solutions $\mathcal{D}^*$ are bi-tropically convex. \edit{Thus, if $\dim \ker (A^{\top}_{E^*})=1$, then $\mathcal{D}^*$ is a polytrope.}
\end{lemma}
\begin{proof}
The feasible region $\mathcal{D}$ is a weighted digraph polyhedron and so by \cref{lem:WDP_are_tropConvex} is bi-tropically convex. Optimal dual solutions $\mathcal{D}^*$ are a face of $\mathcal{D}$ and so by \cref{lem:faces_of_WDP} remain bi-tropically convex. \edit{The final statement concerns boundedness of $\mathcal{D}^*$ after passing to the quotient $\TPT{n}$. The recession cone of $\mathcal{D}^*$ is given by $\ker (A^{\top}_{E*})$, which always contains $\one$. Hence, if $\dim \ker (A^{\top}_{E^*})=1$, then $\ker (A^{\top}_{E*})=\mathrm{span}(\one)$. Since $\mathrm{span}(\one)$ collapses to a point in the quotient, the image of the recession cone in $\TPT{n}$ is trivial. Therefore the image of $\mathcal{D}^*$ contains no nontrivial recession directions and is bounded in $\TPT{n}$, as required.}
\end{proof}

\begin{remark}
    Since polytropes are \edit{tropical} simplices, the Kleene star associated with $\mathcal{D}^*$ provides the complete set of tropical vertices, i.e., generators of $\mathcal{D}^*$.
\end{remark}

\section{Computing Tropical Fermat--Weber Points}\label{sec:Compute_FW}

We first align some notation. We keep $x$ as the flow variable in $\mcf$ and use $\pi^*\in \Pi \subset \TPT{n}$ to denote a tropical Fermat--Weber point. This change is intended to help underscore the duality relation, and emphasize that the space of node potentials has a clean interpretation as a subset of a tropical projective torus. We also introduce index $k\in[m]$ which is an alias for $i$, and is helpful for distinguishing between variables associated with the $\min$ vs. $\max$ components of $f_\nset(\cdot)$.

Recall the functions $\Phi(\cdot)$ and $\Psi(\cdot)$ from the proof of \cref{prop:f_is_convex}. We assign variables $\phi_i=\Phi_i(\cdot)$ for $i\in [m]$ and \edit{$\psi_k=\Psi_k(\cdot)$} for $k\in [m]$ using the classical epigraph reformulation from linear programming. We also introduce the $n\times m$ matrix $Q$ where each column is given by the \emph{negative} of a point in the sample. Since every $p_i\in \nset$ is finite, $-p_i\in \TPT{n}$, so $Q\subset \TPT{n}$ is well-defined. Our new variables introduce constraints according to $\phi_i\geq \pi_j+Q_{ji}$ for $j\in [n]$ and $\psi_k \leq \pi_j+Q_{jk}$ for $j\in[n]$. Finally, define $\delta\mydef (\sum_i\alpha_i-\sum_k\beta_k)/n$, where $\alpha_i$ and $\beta_k$ are as defined in the proof of \cref{prop:f_is_convex}.

\begin{remark}
    Defining $Q$ by taking inverses allows functions $\phi_i$ and $\psi_k$ to be viewed (within their respective semirings) as tropical linear forms, i.e. $q_i\odot \pi$ where $q_i$ is the $i$-th column of $Q$. This provides a cleaner setup for introducing tropical hyperplane arrangements in \cref{sec:subgradient}.
\end{remark}

By negating the objective function and switching to a maximization problem, we arrive at the following primal linear program ($\mathrm{P}$) along with its associated dual ($\mathrm{P'}$). We omit the constant term from the objective function and include the dual variables associated with each block of constraints, as in e.g. $[\mathbf{x}_{ij}]$.
\begingroup
\setlength{\jot}{2pt}
\begin{equation*}
\begin{aligned}
\program{P}\quad 
& \underset{\pi,\,\phi,\,\psi}{\text{minimize}} 
  && \sum_{i=1}^m -\alpha_i \phi_i
     + \sum_{j=1}^n \delta \pi_j 
     + \sum_{k=1}^m \beta{\edit{_k}} \psi_k\\
& \text{subject to} 
  && \pi_j - \phi_i \le -Q_{ji}
  && i\in[m],\;j\in[n] 
  \quad \,[\mathbf{x}_{ij}] \\
& 
  && \psi_k - \pi_j \le Q_{jk}
  && j\in[n],\; k\in[m]
  \quad \edit{[\mathbf{y_{jk}}]} .\\
\program{P'}\quad
& \underset{x\edit{,\,y\,}\ge0}{\text{minimize}} 
  && \sum_{i=1}^m \sum_{j=1}^n
     -Q_{ji}x_{ij} + \sum_{j=1}^n\sum_{k=1}^m Q_{jk}\edit{y_{jk}} \\
& \text{subject to} 
  && \textstyle \sum_{j=1}^n -x_{ij} = -\alpha_i
  && i\in[m]\,
  \quad [\mathbf{\phi}_i] \\
& 
  && \textstyle \sum_{k=1}^m \edit{y_{jk}} - \sum_{i=1}^mx_{ij} = \delta
  && j\in[n] \,\,
  \quad [\mathbf{\pi}_j] \\
& 
  && \textstyle \sum_{j=1}^n \edit{y_{jk}} = \beta_k
  && k\in[m]
  \quad [\mathbf{\psi}_k] .
\end{aligned}
\end{equation*}
\endgroup

We see that ($\mathrm{P}'$) is a min-cost flow problem on $N=2m+n$ nodes and $M=2mn$ arcs. We interpret $[\mathbf{\phi}_i]$ as supply nodes, $[\mathbf{\psi}_k]$ as demand nodes, and $[\mathbf{\pi}_j]$ as transhipment nodes who provide supply or demand as required to balance the system. Arcs $[\mathbf{\phi}_i]\rightarrow [\mathbf{\pi}_j]$ and $[\mathbf{\pi}_j]\rightarrow [\mathbf{\psi}_k]$ exist for all $i\in[m]$, $j\in[n]$, and $k\in [m]$ with costs $(-Q^{\top})_{ij}=-Q_{ji}$ and $Q_{jk}$ respectively. \Cref{fig:mcf2} depicts the $\mcf$ formulation. 
\edit{Note that while we use the additional $y_{jk}$ variables for clarity in the formulations above, the reader should assume that any generic use of $x$ or $x^*$ still refers to \emph{all} flow variables, and not only the left-most half.} 

\edit{It is worth mentioning that the digraph $\mathcal{G}$ representing $(\mathrm{P'})$ has considerably stronger conditions than those provided by our earlier assumptions. In particular, $\mathcal{G}$ is acyclic, which implies \cref{assum:G_bounded} and guarantees feasibility of $(\mathrm{P})$. Additionally, $\mathcal{G}$ is uncapacitated, which, when combined with \cref{assum:G_feasible}, leads to a straightforward proof of feasibility for $(\mathrm{P'})$ using, e.g., a theorem of \cite{Gale1957}.}

\begin{figure}
    \centering
    \begin{minipage}{0.52\textwidth}
        \centering
        \begin{tikzpicture}
        \Vertex[x=0,y=0.8,label={$\phi_m$}]{ym}
        \Vertex[x=0,y=2.2,label={$\phi_1$}]{y1}
        \node[] at (-0.8,3.3) {$[-\alpha_i]$};
        \node at (0,1.6) {$\vdots$};
        \Vertex[x=2,y=0.8,label={$\pi_n$}]{xn}
        \Vertex[x=2,y=2.2,label={$\pi_1$}]{x1}
        \node[] at (2,3.9) {$[\delta]$};
        \node at (2,1.6) {$\vdots$};
        \Vertex[x=4,y=0.8,label={$\psi_m$}]{zm}
        \Vertex[x=4,y=2.2,label={$\psi_1$}]{z1}
        \node[] at (4.8,3.3) {$[\beta_k]$};
        \node at (4,1.6) {$\vdots$};
        \Edge[Direct,bend=10,opacity=1,style={solid}](y1)(x1)
        \Edge[Direct,bend=10,opacity=1,style={solid}](y1)(xn)
        \Edge[Direct,bend=-10,opacity=1,distance=0.6](ym)(x1)
        \Edge[Direct,bend=-10,opacity=1,style={solid}](ym)(xn)
        \Edge[Direct,bend=10,opacity=1,style={solid}](x1)(z1)
        \Edge[Direct,bend=10,opacity=1,style={solid}](x1)(zm)
        \Edge[Direct,bend=-10,opacity=1,distance=0.2](xn)(z1)
        \Edge[Direct,bend=-10,opacity=1,distance=0.2](xn)(zm)
        \Vertex[x=0,y=3.3,label={$\phi_i$}]{uk}
        \Vertex[x=2,y=3.3,label={$\pi_j$}]{xi}
        \Edge[Direct](uk)(xi)
        \node[] at (1,3.65) {$(-Q^{\top})_{ij}$};
        \node[] at (1,3.0) {$x_{ij}$};
        \Vertex[x=4,y=3.3,label={$\psi_k$}]{lk}
        \Edge[Direct](xi)(lk)
        \node[] at (3,3.65) {$Q_{jk}$};
        \node[] at (3,3.0) {\edit{$y_{jk}$}};
    \end{tikzpicture}
    \end{minipage}
    \hfill
    \begin{minipage}{0.46\textwidth}
        \centering
        \[
        \begin{blockarray}{ccccc}
            & \left[\phi\right] & \left[\pi\right] & \left[\psi\right] &\\
        \begin{block}{c(ccc)c}
        \left[-\phi\right] & \mathbb{I} & -Q^{\top} & -Q^{\top}\!\odot\! Q & \alpha\\
        \left[-\pi\right] & \infty & \mathbb{I} & Q & \delta^-\\
        \left[-\psi\right] & \infty & \infty & \mathbb{I} &  0\\
        \end{block}
        & 0 & \delta^+ & \beta &
        \end{blockarray} 
        \]
    \end{minipage}
    \caption{Left: The min-cost flow network for the tropical Fermat--Weber problem. Right: The corresponding optimal transport cost matrix.  $\mathbb{I}$ denotes the tropical identity matrix of appropriate size, $\delta^+=\max(\delta,0)$, and $\delta^-=\max(-\delta,0)$. We can further eliminate the unnecessary $[\psi]$-rows and $[\phi]$-columns which have zero as their corresponding row/column requirements. When $\delta=0$ (as in the symmetric $\fw$ problem), we need consider only the top-right block. Note the negation of the row marginal variables to align our previous variables to a standard optimal transport formulation.}
    \label{fig:mcf2}
\end{figure}

Several strongly polynomial algorithms exist to solve the $\mcf$ problem with arbitrary real-valued costs/capacities, e.g. \citep{tardos_1985_MCC, orlin1993faster, vegh2012strongly}. The Network Simplex algorithm \citep{AMO_NetworkFlows}, although not polynomial-time bounded, also is known to exhibit good practical performance. More recently, interior point methods leveraging sparse approximations and Laplacian solvers with almost-linear-time complexity have been discovered, e.g. \citep{ChenKyngLiuPengGutenbergSachdeva2023}.

With an optimal flow solution $x^*$ to $\mcf$ in hand, we can construct the residual network $\mathcal{G}_{x^*}$ from $\mathcal{G}$ in the standard manner. By \cref{lem:faces_of_WDP} the set of feasible node potentials satisfying the distance inequalities implied by $\mathcal{G}_{x^*}$ form a bi-tropically convex set. We compute the closure of this set by solving the all-pairs shortest path problem over the nodes in $\mathcal{G}_{x^*}$, which yields the Kleene star matrix $D^*$ of optimal node potentials. By the results of \citep{Sergeev2007_closures, JoswigKulas_Trop_Ordinary_Combined}, the rows of $D^*$ are the min-plus tropical vertices of $\mathcal{D}^*$.

The final step requires that we project these tropical vertices from $\TPT{N}$ onto $\TPT{n}$. We make use of the following lemma, referring to its reference for the proof.

\begin{lemma}[\cite{joswigBook}, Lemma 6.6]\label{lem:WDP_proj}
    Let $\mathcal{D}$ be a weighted digraph polyhedron and $D$ be the $r\times r$ matrix with entries in $\Rmin$ associated with $\mathcal{D}$. Additionally, let $D^*$ denote the Kleene star of $D$ and $\mathcal{D}^*$ the weighted digraph polyhedron of $D^*$. Define $\mathrm{proj}_I$ as the map that projects onto coordinates $[r]\setminus I$ for $I\subseteq[r]$. Then the image of $\mathcal{D}^*$ under the linear projection $\mathrm{proj}_I$ is the weighted digraph polyhedron associated to $D^* / I$, where \edit{$D^*/I$ is the matrix arising from $D^*$} by removing the rows and columns which are indexed by elements in $I$.
\end{lemma}

Intuitively, this lemma states that weighted digraph polyhedra are closed under orthogonal projection. In our case, the tropical vertices of the set $\fw(\nset)$ are given as the rows of the $n\times n$ submatrix of $D^*$ obtained by keeping only the rows/columns associated with the $[\mathbf{\pi}_j]$ nodes. Practically, this means we need only compute the shortest paths over $\mathcal{G}_{x^*}$ for every ordered pair $\{(j,j')\mid j, j' \in [n]\}$ rather than the full matrix $D^*$. Computing these shortest path distances yields $D^*$, whose rows are the tropical vertices of $\fw(\nset)$.

The following result follows immediately.
\begin{theorem}
    For any finite sample $\nset \subset \TPT{n}$, every weighted tropical Fermat--Weber point is isomorphic to an optimal node potential vector of the corresponding min-cost flow problem, and the set of all such tropical Fermat--Weber points forms a polytrope in $\TPT{n}$.
\end{theorem}
\begin{proof}
    Any fixed $\pi\in \TPT{n}$ yields a feasible solution in $(\mathrm{P})$ through appropriate choices of $\phi_i$ and $\psi_k$, so $(\mathrm{P})$ is feasible. Since $f_\nset(\cdot)$ is coercive, the set of optimal $\pi$ is bounded, and the feasible region of $(\mathrm{P})$ imposes corresponding bounds on $\phi_i$ and $\psi_k$. Hence, the set of optimal solutions of $(\mathrm{P})$ is bounded. The feasible region of $(\mathrm{P})$ is a weighted digraph polyhedron and \cref{lem:MCF_duals_are_tropConvex} shows that its optimal face forms a bi-tropically convex set. \cref{lem:WDP_proj} ensures that the orthogonal projection to $\TPT{n}$ remains bi-tropically convex. Since projection to a lower-dimensional subspace preserves boundedness, we have a bounded, bi-tropically convex set, i.e., a polytrope. By standard linear programming duality, optimal solutions of $(\mathrm{P})$ admit an interpretation as node potential vectors for the associated dual min-cost flow problem.
\end{proof}

Since $\mathcal{G}$ (the graph associated to our $\mcf$ problem) is acyclic and uncapacitated, our dual program has an equivalent formulation as a transportation problem, see e.g., \citep{AMO_NetworkFlows}. The association between tropical Fermat--Weber problems and the transportation problem was first identified by \cite{Com_neci_2023} in the special cases $\gamma=\one$ and $\gamma=\boldsymbol{0}$. In these instances, $\mathcal{G}$ arises as an optimal transport problem directly. In the general case, transformation into the optimal transport formulation is straightforward.

The transformation requires that we compute $C=(c_{ij}) \in \R^{N\times N}$, the matrix of shortest path distances between all pairs of nodes in $\mathcal{G}$. Since $\mathcal{G}$ is acyclic, negative cycles can't exist, implying the Kleene star is well-defined. Since $\mathrm{diam}(\mathcal{G})\leq 2$, the only non-trivial pairwise distances we need to compute are those between nodes $\phi_i\rightarrow\psi_k$ for $i,k\in[m]$. In tropical terms, $C_{ik}=(-Q^{\top}\odot Q)_{ik}$, which can be computed in $\mathcal{O}(m^3)$. The complete cost matrix used for the optimal transport $\ot$ formulation is shown on the right in \cref{fig:mcf2}.


\begin{remark}
    When $\gamma=\tfrac{1}{2}\one, \,w=\one$, the cost matrix simplifies to an $m \times m$ assignment problem, for which the Hungarian algorithm can provide a solution in $\mathcal{O}(m^3)$. For rational weights, this also holds under the more general setting after suitably expanding the problem via scaling and node replication. Thus, every tropical Fermat--Weber problem with rational weights has a formulation as an assignment problem.
\end{remark}

Conversion from $\mcf$ to $\ot$ increases the density of the underlying graph (by adding arcs $[\phi]\rightarrow [\psi]$). This suggests that even if ignoring the $\mathcal{O}(m^3)$ cost to form $C$, there is little computational benefit to solving $\ot$ rather than $\mcf$ directly. Nonetheless, connections to optimal transport and the assignment problem  provide insights on the geometry of the tropical Fermat--Weber problem. For example, the \emph{tropical volume} of a square matrix $A\in \R^{n\times n}$ is defined by solving the assignment problem on $A$, see \citep{joswigBook}. Indeed, since volume generalizes the notion of distance to higher dimensions, its connection to a location problem is natural. 

\begin{corollary}
    The tropical volume $\mathrm{tvol}$ of the tropical Fermat--Weber set $\mathrm{tvol}(\fw(\nset))$ is equal to 
    \[
    \mathrm{tvol}(\fw(\nset))=\bigoplus_{\sigma \in (\mathrm{Sym(n)-\sigma^*})}\sum D^*_{i,\sigma(i)} - \bigoplus_{\sigma \in \mathrm{Sym}(d)}\sum D^*_{i,\sigma(i)}
    \]
    where $D^*$ is the $n\times n$ Kleene star matrix whose rows are the min-plus tropical extreme points of $\fw(\nset)$, and $\sigma^*$ is an optimal permutation of the assignment problem on $D^*$.
\end{corollary}

In the statistical physics literature, this volume is interpreted as an energy gap, which \cite{YuilleKosowsky1994InvisibleHand} use in analyzing the convergence of their algorithm for solving the linear assignment problem over integer matrices. In particular, when this energy gap is known in advance, they show that it can be used to inform a sufficiently low temperature for a steepest descent scheme minimizing a ``\emph{P}-energy.'' By performing a fixed number of descent steps at this temperature, they obtained a required accuracy threshold, enabling round-off to the nearest integer which is optimal. In general, however, this energy gap is unknown apriori and so other approaches are required.

\section{Iterative Methods}\label{sec:subgradient}

In addition to combinatorial and simplex-style algorithms, convex optimization problem admits, under standard regularity conditions, iterative first-order methods that converge to a global optimum or generate $\epsilon$-optimum solutions. Such iterative schemes can be useful in certain settings, for example, when any $\epsilon$-approximate solution is acceptable, or if operating under certain computational constraints. These include subgradient methods, proximal methods, and related first-order techniques \citep{bertsekas2015convex}. We formulate the tropical Fermat--Weber subdifferential of the objective as a function of $\pi\in \TPT{n}$ and describe it from the tropical perspective.

For a linear program $\min_\pi \{b^{\top}\pi \mid A^{\top}\pi\leq c\}$, consider its indicator function formulation $F(\pi)=b^{\top}\pi+\mathbbm{1}_\mathcal{D}(\pi)$, where $\mathcal{D}=\{\pi\mid A^{\top}\pi\leq c\}$ is the \edit{polyhedron} defining the feasible region and $\mathbbm{1}_\mathcal{D}$ denotes the (tropical) indicator function. That is, $\mathbbm{1}_\mathcal{D}(\pi)=0$ if $\pi\in \mathcal{D}$ and $\mathbbm{1}_\mathcal{D}(\pi)=\infty$ otherwise. The \emph{subdifferential} is $\partial F(\pi)=b+N_\mathcal{D}(\pi)$ where $N_\mathcal{D}(\pi)$ is the normal cone to $\mathcal{D}$ at $\pi$. For feasible $\pi$, $N_\mathcal{D}(\pi)=\{Ax \mid x\geq0, x_i(a_i\pi-c_i)=0\}$, where $a_i$ is the $i$-th column of $A$, and so $\partial F(\pi)=b+Ax$ for active multipliers $x$.

Recall that for the tropical Fermat--Weber linear program $(\mathrm{P})$, any $\pi\in \TPT{n}$ corresponds to a feasible solution once we fix $\phi_i=\max_j(\pi_j+Q_{ij})$ for $i\in [m]$ and $\psi_k=\min_j(\pi_j+Q_{kj})$ for $k\in [m]$. Doing so immediately identifies which constraints in $(\mathrm{P})$ are active; equivalently, the face of $\mathcal{D}$ corresponding to $\pi$, which we denote as $\mathcal{D}(\pi)$. Remove from $\mathcal{D}(\pi)$ all inactive constraints and construct the corresponding equality subgraph $\mathcal{G}|_\pi$, which we say is induced by $\pi$. Note that $\mathcal{G}|_\pi$ is equivalent to the graph $\mathcal{G}$ of $\mcf$ that enforces flow $x_{ij}=0$ for any arc corresponding to an inactive constraint in $\mathcal{D}(\pi)$. Similar to our notation in \cref{sec:basics}, we use $E^a\subseteq E$ to denote the arc set of $\mathcal{G}|_\pi$. For any flow $x$ on $\mathcal{G}|_\pi$, the degree to which flow conservation across nodes is violated is precisely the subdifferential at $\pi$. Since each node in $[\mathbf{\phi}]$ has at least one outgoing arc, and each node in $[\mathbf{\psi}]$ has at least one incoming arc, we can always select flows that satisfy node balance for our auxiliary variables. For example, if we assign flow $x_{ij}=\alpha_i$ for $j$ the lexicographically first outgoing arc from node $\phi_i$ (in $\mathcal{G}|_\pi$) and assign $x_{\cdot j}=0$ to all others, constraints associated with $[\phi_i]$ in $(\mathrm{P'})$ are satisfied. Therefore, there always exists a subgradient $g(\pi)\in \partial F(\pi)$ for which the only non-zero elements correspond to nodes in $[\mathbf{\pi}]$. Selecting a subgradient in this manner amounts to a projected subgradient method in which the projection is implicit based on the assignment of dual (flow) variables. This basic approach yields \cref{alg:basic_fw_gd}.

\begin{algorithm}
\caption{Tropical Fermat--Weber Subgradient Descent}\label{alg:basic_fw_gd}
\begin{algorithmic}
    \State \textbf{Input}: Sample $\nset = \{\sample \}\subset\TPT{n}$, weights $\alpha,\,\beta \in \R^m_{\geq0}$, max iterations $K$, initial point $\pi^0\in \TPT{n}$, and step size schedule $\eta=(\eta^1,\ldots,\eta^K)$.
    \State \textbf{Output}: $\pi^k\in \fw(\nset)$ or final iterate $\pi^K$.
    \For{$k=1,\ldots, K$}
        \State Set $x\leftarrow \boldsymbol{0}$.
        \State Identify the set of active constraints $\mathcal{D}(\pi^k)$.
        \For{$i=1,\ldots, m$} assign flow so that $\sum_j\bigl(x_{ij}\mid (i,j)\in E^a\bigr)=\alpha_i$.
        \EndFor
        \For{$k=1,\ldots, m$} assign flow so that $\sum_j\bigl(\edit{y_{jk}}\mid (j,k) \in E^a\bigr)=\beta_k$.
        \EndFor
        \State Compute $g(\pi^k)=c-Ax$.
        \If{$g(\pi^k)=\boldsymbol{0}$} stop. \EndIf
        \State $\pi^{k+1} \leftarrow \pi^k - \eta^kg(\pi^k)$ .
    \EndFor
    \State \textbf{return} $\pi^k$
\end{algorithmic}
\end{algorithm}

The subgradients we consider here are quite similar to those in \citep{TalbutMonod2025}, in which the authors propose a more general \emph{tropical descent} framework for problems over $\TPT{n}$, where they include tropical location problems as one such problem class. In that work, the authors focus on tropically convex problems which are, in general, not classically convex. In addition, their Algorithm 1 imposes binary descent directions in a way that is conceptually more similar to the methods used for $L^{\sharp}$-convex minimization (e.g., \citep{Murota2022}) or non-divisible product-mix auctions (e.g., \citep{BaldwinGoldbergKlempererLock2023}). Since our problem is convex, we leverage the more general setting provided by classical convex analysis and impose no such restrictions on the (sub)gradients.

With an appropriate step size schedule, \cref{alg:basic_fw_gd} will converge to some $\pi^*\in \fw(\nset)$. For example, the classical diminishing step-size schedule $\eta_k=2/(k+2)$ yields $\mathcal{O}(1/k)$ convergence under standard convexity assumptions, see \citep{bertsekas2015convex}. It is well-known, however, that convergence under such schedules can be extremely slow. Furthermore, we have not specified how one should assign flows in cases where the maxima/minima used to fix $\phi_i$ and $\psi_k$ are attained over more than one index $j\in[n]$. To help inform these considerations, we consider the problem's underlying geometry, which turns out to be quite intuitive when viewed from the tropical perspective. 

\textbf{Tropical Geometry Part II.} For $v\in \Rmin^n$, the \emph{support} $\mathrm{supp}(v)$ is the set of indices $j \in [n]$ for which $v_j$ is finite. For $v\in \Rmin^n$ with $|\mathrm{supp}(v)|\geq2$, the tropical linear form $v_1\odot \pi_1 \oplus \cdots \oplus v_n \odot \pi_n$, defines a \emph{tropical hyperplane} $\mathcal{H}_v$ consisting of all points $\pi \in \TPT{n}$ for which \edit{$\min_i(\pi_i+v_i)=\min (\pi_j+v_j)$ for some $i\neq j$}. The requirement on the support of $v$ ensures that the hyperplane is non-empty. When $v$ has full support, i.e. $v\in \TPT{n}$, then there is a point $\pi=-v+\R\one$, unique up to scalar translation, for which all terms $\pi_j+v_j$ are minimized simultaneously. We call this point the \emph{apex} of $\mathcal{H}_v$. 

Tropical hyperplanes partition $\TPT{n}$ into \emph{sectors}. For $v\in \Rmin^n$ and $j\in \mathrm{supp}(v)$, the $j$-th (closed) sector (of the partition of $\TPT{n}$ by $\mathcal{H}_v$) is the set 
\[
S_v^j\mydef \bigcap_{i\in \mathrm{supp}(v)\setminus j}\{\pi\in \TPT{n} \mid \pi_j+{\edit{v_j}} \leq \pi_i+{\edit{v_i}}\} .
\]
Enforcing strict inequality in the definition above defines the $j$-th \emph{open} sector, $S_v^{\circ,j}$. For $j\in \mathrm{supp}(v)$ it is clear that $\pi\in S_v^j \iff j\in \argmin_i(\pi_i+v_i)$.

The partition of $\TPT{n}$ given by the inequalities implied by $S_v^j$ again defines a weighted digraph polyhedron, and so we see that tropical sectors are bi-tropically convex. Recalling the constraints in $(\mathrm{P})$, we see 
\[
\psi_i=\min_j(\pi_j-p^i_j)=(Q_{1i} \odot \pi_1 \oplus \cdots \oplus Q_{ni} \odot \pi_n)\leq \pi_j+Q_{ji},  \quad j\in [n]
\]
is simply the tropical hyperplane $\mathcal{H}_{q_i}$ with apex given by $-q_i=p_i$. By considering the full set of points $p_i\in \nset$, we get a tropical hyperplane arrangement $\mathcal{T}(\nset)$ (union of several hyperplanes) whose common refinement is a polyhedral subdivision of $\TPT{n}$. Each cell $\mathcal{C}$ of this subdivision is an intersection of tropical sectors (and so is bi-tropically convex) which encodes the combinatorial structure of the $\mathcal{T}(\nset)$. That is, for any point $\pi \in \mathcal{C}$ in the cell, we take the $n$-tuple whose $j$-th element records the indices for points $p_i\in \nset$ whose $j$-th closed sector contains $\pi$. We call this tuple the \emph{covector} of $\mathcal{C}$. The set of all such covectors yields the \emph{covector decomposition}, $\CovDec(\nset)$, generated by $\mathcal{T}(\nset)$, which encodes the combinatorial structure of the tropical hyperplane arrangement. \cite{Ardila} showed that this subdivision is combinatorially isomorphic to the complex of interior faces of the regular subdivision of a product of simplices.
\begin{remark}
    For an unweighted asymmetric $\fw$ problem, the set $\fw(\nset)$ is precisely the cell of $\CovDec(\nset)$ corresponding to the set of \emph{tropical centerpoints} (of $\nset$). Assuming $m\geq n$, any tropical hyperplane with apex $-\pi$ for $\pi\in\fw(\nset)$ provides a partition of $\nset$ according to the hyperplane's $n$ sectors. For integer weighted sample $\nset$, points can be split into multiple lower-weighted copies (e.g. via 2-adic decomposition) to form a larger sample $\nset^+$ such that $x\in\fw(\nset)$ is a \emph{tropical $n$-Tverberg point} of $\nset^+$. This is precisely the method outlined in \citep{TokuyamaNakano1995} for the Hitchcock transportation problem. For definitions of tropical centerpoints and Tverberg points, see \citep{JaggiKatzWagner2010}. 
\end{remark}

For any $\pi\in\TPT{n}$, the covector of the maximal cell containing $\pi$ encodes the active constraints of \edit{$[\mathbf{y}_{jk}]$} in $(\mathrm{P})$. By taking max-plus analogues, we see that the max-plus covector of the maximal cell containing $\pi$ encodes the active constraints of $[\mathbf{x}_{ij}]$. Thus, the subgradient descent methodology we outlined earlier essentially corresponds to tracking and updating min-plus and max-plus covectors over a common refinement of $\TPT{n}$ given by the union of the min-plus tropical hyperplane arrangement $\mathcal{T}_{\min}(\nset)$ and the max-plus tropical hyperplane arrangement $\mathcal{T}_{\max}(\nset)$. We call this joint subdivision the \emph{bi-tropical covector decomposition} of $\nset$, denoted as $\BCD(\nset)$. This interpretation also leads to the following description of $\fw(\nset)$.

\begin{theorem}
    For any finite sample $\nset\subset \TPT{n}$, the polytrope corresponding to the set of tropical Fermat--Weber points is given precisely by a cell of the bi-tropical covector decomposition $\BCD(\nset)$. 
\end{theorem}
\begin{proof}
    Since $(\mathrm{P})$ feasible and bounded, the set of optimal solutions forms a unique face of the feasible region given explicitly as the set of constraints active at $\pi^*\in \fw(\nset)$. Each active constraint in this set corresponds bijectively with the tropical sector indicated by the non-zero indices of the constraint. Applying this to the full set of active constraints equates to intersecting the corresponding set of tropical sectors. This is precisely the refinement given by $\BCD(\nset)$. 
\end{proof}

The construction of $\BCD(\nset)$ is closely related to the \emph{ambitropical cones} and \emph{ambitropical hulls} of \citep{AkianGaubertVannucci2023}. 
\begin{definition}[Ambitropical Cones and Hulls, \cite{AkianGaubertVannucci2023}]
    An \emph{ambitropical cone} is a non-empty additive cone $C$ of $\R^n$ such that $C$ is a lattice in the induced order of $(\R^n,\leq)$. For $U\subset \R^n$, $\tilde{U}\subset \R^n$ is an \emph{ambitropical hull} of $U$ if it is a closed ambitropical cone which is a superset of $U$ and is minimal with respect to inclusion.
\end{definition}
The min-plus and max-plus convexity of a bi-tropically convex set together imply the lattice structure required of an ambitropical cone, so in our setting these two notions are equivalent. This implies that each non-empty closed cell of $\BCD(\nset)$ is an ambitropical cone, and hence, by the results of \citep{AkianGaubertVannucci2023}, also a \emph{hyperconvex set} (see \citep{Isbell1964,Dress1984}) and also a fixed-point set for a particular \emph{Shapley operator}. Since (abstract) Shapley operators are dynamic programming operators used to compute the value function of a zero-sum game, our $\BCD$ construction demonstrates the linkage between a tropical Fermat--Weber problem and that of finding a price vector $\pi$ that achieves equilibrium in a particular competitive market. Indeed, this analogy is precisely the one given by the Auction Algorithm of \citep{BertsekasEckstein1988} which is used to solve the Assignment Problem.

The following example illustrates these bi-tropical objects. 
\begin{example}\label{ex:example}
 Let $\nset=\{p_1,p_2,p_3,p_4\}$ be the four points in $\TPT{3}$ given by $p_1=(0,0,5)^{\top}$, $p_2=(0,1,2)^{\top}$, $p_3=(0,3,0)^{\top}$, and $p_4=(0,3,6)^{\top}$ as shown in  \cref{fig:Example}. For simplicity we consider the symmetric, unweighted case. That is, \edit{$\gamma_i=1/2$} and $w_i=1$ for $i=1,\ldots,4$. For $\pi=(0,3,2)^{\top}$ we show the equality subgraph $\mathcal{G}|_\pi$ that defines the subdifferential $\partial F(x)$. $\fw(\nset)$ is the gray shaded cell. The tropical vertices are marked by white circles at $p_2$, $(0,2,4)^{\top}$, and $(0,3,5)^{\top}$.

 \begin{figure}[H]
    \centering
    \begin{tikzpicture}[node distance={30mm}, thick, main/.style = {draw, circle}]
        \fill[{gray!20}] (1,2) -- (3,4) -- (3,5) -- (2,5) -- (1,4);
        \node at (0,5)[circle,fill,inner sep=1.5pt]{};
        \node at (1,2)[circle,fill,inner sep=1.5pt]{};
        \node at (3,0)[circle,fill,inner sep=1.5pt]{};
        \node at (3,6)[circle,fill,inner sep=1.5pt]{};
        \draw [dashed,->](0,5) to (0,7);
        \draw [dashed,->](0,5) to (4.7,5);
        \draw [dashed,->](0,5) to (-2,3);
        \draw [->](0,5) to (0,-1);
        \draw [->](0,5) to (-2,5);
        \draw [->](0,5) to (2,7);
        \draw [dashed,->](1,2) to (1,7);
        \draw [dashed,-](1,2) to (3,2);
        \draw [dashed,->](1,2) to (-2,-1);
        \draw [->](1,2) to (1,-1);
        \draw [->](1,2) to (-2,2);
        \draw [->](1,2) to (6,7);
        \draw [dashed,->](3,0) to (4,0);
        \draw [dashed,->](3,0) to (2,-1);
        \draw [->](3,0) to (-2,0);
        \draw [->](3,0) to (3,-1);
        \draw [->](3,0) to (3.8,0.8);
        \draw [dotted](3,0) to (3,1);
        \draw [dashed,->](3,6) to (3,7);
        \draw [dashed,->](3,6) to (6.5,6);
        \draw [dashed,->](3,6) to (-2,1);
        \draw [->](3,6) to (-2,6);
        \draw [dotted](3,6) to (3,2);
        \draw [->](3,6) to (4,7);

        \node[]  at (0.3,4.7) {$p_1$};
        \node[]  at (1.3,1.7) {$p_2$};
        \node[]  at (3.3,-0.3) {$p_3$};
        \node[]  at (3.3,5.7) {$p_4$};
        
        \node[] at (2.8,2.2) {$\pi$};
        \draw [draw=red,->,dotted](3,2) to (3,2-3/3);
        \draw [draw=red,->,dotted](3,2) to (3+5/3,2+1/3);
        \draw [draw=red,->](3,2) to (3+3/3,2);
        \draw [draw=red,->](3,2) to (3+1/3,2-1/3);
        \node[]  at (3.7,1.5) {$g(\pi)_{\mathrm{opt}}$};
        \node[]  at (4.55,1.9) {$g(\pi)_{\mathrm{lex}}$};
        
        \node at (3,2)[circle,draw=black,fill=white,inner sep=1.5pt]{};
        \node at (3,5)[circle,draw=black,fill=white,inner sep=1.5pt]{};
        \node at (1,4)[circle,draw=black,fill=white,inner sep=1.5pt]{};
        \node at (1,2)[circle,draw=black,fill=white,inner sep=1.5pt]{};
        
        \begin{scope}[xshift=6cm, yshift=3cm]
            \draw [](0,0.00) to (1.5,2);
            \draw [](0,0.00) to (1.5,1);
            \draw [](0,0.66) to (1.5,0);
            \draw [](0,1.33) to (1.5,1);
            \draw [](0,2.00) to (1.5,1);
            
            \draw [solid](1.5,0) to (3,0.00);
            \draw [solid](1.5,1) to (3,0.66);
            \draw [solid](1.5,2) to (3,0.66);
            \draw [solid](1.5,0) to (3,1.33);
            \draw [solid](1.5,2) to (3,1.33);
            \draw [solid](1.5,0) to (3,2.00);

            \node at (0, 0.00)[circle,fill,inner sep=1.5pt]{};
            \node at (0, 0.66)[circle,fill,inner sep=1.5pt]{};
            \node at (0, 1.33)[circle,fill,inner sep=1.5pt]{};
            \node at (0, 2.00)[circle,fill,inner sep=1.5pt]{};
            
            \node at (1.5, 0)[circle,draw=black,fill=white,inner sep=1.5pt]{};
            \node at (1.5, 1)[circle,draw=black,fill=white,inner sep=1.5pt]{};
            \node at (1.5, 2)[circle,draw=black,fill=white,inner sep=1.5pt]{};
            
            \node at (3, 0.00)[circle,fill,inner sep=1.5pt]{};
            \node at (3, 0.66)[circle,fill,inner sep=1.5pt]{};
            \node at (3, 1.33)[circle,fill,inner sep=1.5pt]{};
            \node at (3, 2.00)[circle,fill,inner sep=1.5pt]{};
            
            \node[] at (-0.4, 0.00) {$\phi_4$};
            \node[] at (-0.4, 0.66) {$\phi_3$};
            \node[] at (-0.4, 1.33) {$\phi_2$};
            \node[] at (-0.4, 2.00) {$\phi_1$};
            \node[] at (3.4, 0.00) {$\psi_4$};
            \node[] at (3.4, 0.66) {$\psi_3$};
            \node[] at (3.4, 1.33) {$\psi_2$};
            \node[] at (3.4, 2.00) {$\psi_1$};
            \node[] at (1.5, 2.5) {$\mathcal{G}|_\pi$};
        \end{scope}
        \begin{scope}[xshift=6cm, yshift=-0.5cm]
            \draw [draw=gray,dotted](0,2.00) to (1.5,2);
            \draw [draw=gray,dotted](0,2.00) to (1.5,1);
            \draw [draw=gray,dotted](0,1.33) to (1.5,1);
            \draw [draw=gray,dotted](0,1.33) to (1.5,0);
            \draw [draw=gray,dotted](0,0.66) to (1.5,2);
            \draw [draw=gray,dotted](0,0.66) to (1.5,0);
            \draw [draw=gray,dotted](0,0.00) to (1.5,2);
            \draw [draw=gray,dotted](0,0.00) to (1.5,1);

            \draw [draw=gray,dotted](3,2.00) to (1.5,2);
            \draw [draw=gray,dotted](3,2.00) to (1.5,0);
            \draw [draw=gray,dotted](3,1.33) to (1.5,2);
            \draw [draw=gray,dotted](3,1.33) to (1.5,1);
            \draw [draw=gray,dotted](3,0.66) to (1.5,2);
            \draw [draw=gray,dotted](3,0.66) to (1.5,1);
            \draw [draw=gray,dotted](3,0.00) to (1.5,1);
            \draw [draw=gray,dotted](3,0.00) to (1.5,0);
            
            \draw [](0,0.00) to (1.5,0);
            \draw [](0,0.66) to (1.5,1);
            \draw [](0,1.33) to (1.5,2);
            \draw [](0,2.00) to (1.5,0);
            \draw [solid](1.5,2) to (3,0.00);
            \draw [solid](1.5,1) to (3,2.00);
            \draw [solid](1.5,0) to (3,1.33);
            \draw [solid](1.5,0) to (3,0.66);

            \node at (0, 0.00)[circle,fill,inner sep=1.5pt]{};
            \node at (0, 0.66)[circle,fill,inner sep=1.5pt]{};
            \node at (0, 1.33)[circle,fill,inner sep=1.5pt]{};
            \node at (0, 2.00)[circle,fill,inner sep=1.5pt]{};
            
            \node at (1.5, 0)[circle,draw=black,fill=white,inner sep=1.5pt]{};
            \node at (1.5, 1)[circle,draw=black,fill=white,inner sep=1.5pt]{};
            \node at (1.5, 2)[circle,draw=black,fill=white,inner sep=1.5pt]{};
            
            \node at (3, 0.00)[circle,fill,inner sep=1.5pt]{};
            \node at (3, 0.66)[circle,fill,inner sep=1.5pt]{};
            \node at (3, 1.33)[circle,fill,inner sep=1.5pt]{};
            \node at (3, 2.00)[circle,fill,inner sep=1.5pt]{};
            
            \node[] at (-0.4, 0.00) {$\phi_4$};
            \node[] at (-0.4, 0.66) {$\phi_3$};
            \node[] at (-0.4, 1.33) {$\phi_2$};
            \node[] at (-0.4, 2.00) {$\phi_1$};
            \node[] at (3.4, 0.00) {$\psi_4$};
            \node[] at (3.4, 0.66) {$\psi_3$};
            \node[] at (3.4, 1.33) {$\psi_2$};
            \node[] at (3.4, 2.00) {$\psi_1$};
            \node[] at (1.5, 2.5) {$\mathcal{G}_{\pi^*}$};
            \node[] at (0.75, 2.3) {$x^*_{ij}$};
            \node[] at (2.25, 2.3) {$x^*_{jk}$};
        \end{scope}
        \node[] at (2,4) {$\Pi^*$};
    \end{tikzpicture}
    \caption{The \edit{bi-tropical} covector decomposition $\BCD(\nset)$ of four points in $\TPT{3}$. Max-plus hyperplanes (solid) and min-plus hyperplanes (dashed) together subdivide the plane. The segment from $p_3$ to $p_4$ is dotted to show it corresponds to the overlap of max-plus (from $p_4$) and min-plus (from $p_3$). That is, points in $\nset$ are not in {\em general position}. The cell corresponding to $\fw(\nset)$ is shaded gray and its three tropical vertices are marked with white circles. The equality subgraph at $\pi$ on the right shows the possible assignments of flow that can be used to compute a subgradient at $\pi$. The normal cone defining the subdifferential $\partial F(\pi)$ is depicted by its extreme rays (dotted red). Two possible subgradients are also shown corresponding to a lexicographic selection rule, and the optimal selection obtained by solving a quadratic minimization problem. The bottom right graph shows the optimal flow solution to $\mcf$, with active flows drawn solid.}
    \label{fig:Example}
\end{figure}
\end{example}

Any $n$-dimensional cell $\mathcal{C}$ of $\BCD(\nset)$ has points in its relative interior that are not contained in any hyperplane of $\mathcal{T}_{\min}(\nset)$ or $\mathcal{T}_{\max}(\nset)$. Such cells are sometimes called \emph{topes}. Equivalently, the min-plus and max-plus sectors containing $\mathcal{C}$ are uniquely defined for each $p_i\in \nset$, and the elements of their corresponding covectors are all singletons. Hence, when $\pi\in \mathrm{relint}(\mathcal{C})$ for a tope $\mathcal{C}$, our rule for subgradient selection is uniquely determined, and any check $\boldsymbol{0}\in \partial F(\pi)$ of optimality conditions proceeds directly. When $\pi$ is contained in a hyperplane, however, selection becomes non-trivial. To help address this, we consider the following two observations. 

First, let us assume that all points in the sample $\nset$ are in {\em general position}, that is, the matrix $Q$ has no $r\times r$-submatrix which is tropically singular for $2 \leq r \leq n$ (recall we assume $m\geq n)$.
\edit{By \emph{tropically singular}, we mean that all of the row (or column) vectors of the corresponding square matrix lie on a tropical hyperplane; see \citep[Proposition 3.4]{joswigBook}.}
When $\nset$ is in general position, no point $p_i$ is contained in hyperplane arrangements given by $\nset \setminus p_i$. That is, $p_i \notin \mathcal{T}_{\min}(\nset \setminus p_i)\cup \mathcal{T}_{\max}(\nset \setminus p_i)$. Thus, by considering the set $\nset\setminus p_i$ rather than $\nset$, the subgradient at $\pi=p_i$ becomes unique per our selection rule. Furthermore, if $\pi=p_i$ is a tropical Fermat--Weber point of $\nset \setminus p_i$, then by \citep[Lemma 8]{Lin2016TropicalFP} $\fw(\nset)=p_i$ is the unique tropical Fermat--Weber point of the original sample. 

\begin{remark}
    In \citep{Lin2016TropicalFP}, the authors only consider the unweighted symmetric tropical Fermat--Weber problem. Since $\drtr$ is non-negative, we see their result holds in the more general context. Intuitively, adding a new point $u$ to $\nset$ will not increase $f(\pi)$ if and only if $u\in \fw(\nset)$. Doing so results in $\fw(\nset\cup u)=u$. 
\end{remark}

Even if $\nset$ is not in general position, any point $p_i$ for which $p_i \notin \mathcal{T}_{\min}(\nset \setminus p_i)\cup \mathcal{T}_{\max}(\nset \setminus p_i)$ can leverage this result. In some instances, pre-computing subdifferentials at points $p_i\in \nset$ can identify $p_i\in \fw(\nset)$ much faster than solving $\mcf$ in its entirety.

For the general situation in which $\pi$ is contained in one or more tropical hyperplanes, we leverage a known correspondence between min-cost flow and max-flow problems. Take $\nset = \{\sample \}\subset\TPT{n}$, $\pi \in \TPT{n}$, and $\mathcal{G}|_\pi$ the equality subgraph of $\mathcal{G}$ at $\pi$. To $\mathcal{G}|_\pi$, add a special ``start'' node $s$ and arcs $(s,\phi_i)$ with cost $c(s,\phi_i)=0$ and capacity $\mu(s,\phi_i)=\alpha_i$ for all $i\in [m]$. Similarly, add a ``terminal'' node $t$ and arcs $(\psi_k,t)$ with cost $c(\psi_k,t)=0$ and capacity $\mu(\psi_k,t)=\beta_k$ for all $k\in [m]$. Ignoring all node supplies/demands, find the maximum flow from $s\rightarrow t$ in this modified equality graph. This problem, $\mathrm{MF}$, always has an optimal solution, and we record the objective function value as $\rho^*$. The following proceeds directly from well-known results in network optimization, see e.g. \citep[Theorem 9.7]{AMO_NetworkFlows}. It's proof is omitted.

\begin{proposition}
     For $\mathrm{MF}$, $\pi\in \fw(\nset)$ if, and only if, $\rho^*=\sum_i\alpha_i$.
\end{proposition}
\begin{remark}
    Since the max-flow instance $\mathrm{MF}$ is an acyclic graph in which every $s\rightarrow t$ path has length at most $4$, the shortest-path distance from $s$ to $t$ is the residual graph is bounded by a constant. In Dinic's algorithm \citep{dinic1970algorithm}, this distance increases after each blocking-flow phase, implying $\mathcal{O}(1)$ phases are required. As each phase can be implemented in $\mathcal{O}(|E|)$ time, the total running time is $\mathcal{O}(|E|)$. Recall that typically $|E|\ll M=2mn$ for most equality graphs $\mathcal{G}|_x$.
\end{remark}

Even when $\pi$ is not contained in a hyperplane, we might want to consider which hyperplane sectors are within some $\theta$-ball of $\pi$, $\mathcal{B}_\theta(\pi)$. Computing the set of such sectors is straightforward. For any point $p$, the set of min-plus sectors within an $\theta$-ball of $\pi$ is computed as $\{j\mid x_j - p_j - \min_k(x_k-p_k)\leq \theta\}$. For max-plus sectors, we similarly have $\{j\mid x_j - p_j - \max_k(x_k-p_k) \geq -\theta\}$. By considering the set of sectors reachable within some $\theta$-ball (in the tropical norm) of $\pi$, we can run a relaxed version of $\mathrm{MF}$ to check a necessary condition for optimality. While this relaxed version clearly isn't sufficient (consider $\theta\rightarrow \infty$), it does provide insight on the suitability of a given step size. In fact, if one knew the exact minimal $\theta$ for which $\fw(\nset) \cap \mathcal{B}_\theta(\pi)\neq \emptyset$, we could solve $\mcf$ on the subgraph containing only those arcs corresponding to hyperplane sectors with a non-empty intersection with $\mathcal{B}_\theta(\pi)$. One might consider combining such checks in conjunction with a line search routine or some other adaptive step-size method, or in eliminating arcs from $\mcf$ which, by the requirement that $\fw(\nset) \subseteq \mathrm{tconv}(\nset)$, need never be considered. For starters, any $\pi$ must have $\dtr(\pi,\fw(\nset))\leq \max_i\{\dtr(\pi,p_i)\}$, providing a conservative initial value for $\theta$.

The basic procedure given in \cref{alg:basic_fw_gd} can be adapted to check optimality criteria and better select subgradients when $\pi$ lies on a tropical hyperplane of $\BCD(\nset)$. On any iteration, after determining the set of active constraints, one can perform a series of checks. If $\pi=p_i$, remove $p_i$ from $\nset$ before continuing. If $\pi$ doesn't lie on any tropical hyperplane, compute the subgradient as usual. Otherwise, solve $\mathrm{MF}$ as discussed previously, and check the optimality criteria. Assuming the criteria fails to be satisfied, select a subgradient according to some rule, for example, by minimizing some norm on the $[\pi]$ node imbalances. The overall procedures is given in \cref{alg:improved_fw_gd}.

\begin{algorithm}
\caption{Tropical Fermat--Weber Subgradient Selection}\label{alg:improved_fw_gd}
\begin{algorithmic}
    \State \textbf{Input}: Sample $\nset=\{\sample \}\subset \TPT{n}$, weights $\alpha,\,\beta \in \R^m_{\geq0}$, and $\pi\in \TPT{n}$.
    \State \textbf{Output}: $g(\pi) \in \partial F_\nset(\pi)$.
    \State Generate the equality subgraph $\mathcal{G}|_\pi$ by deleting inactive arcs.
    \State Delete any $[\mathbf{\phi}]$ or $[\mathbf{\psi}]$ nodes of degree $n$.
    \If{all nodes in $\{[\mathbf{\phi}]\uplus [\mathbf{\psi}]\}$ have degree $1$} assign flows $x$ as in \cref{alg:basic_fw_gd}.
        \State Compute $g(\pi)=c-Ax$.
    \Else 
        \State Solve $\mathrm{MF}$ over $\mathcal{G}|_x$ to obtain maximum flow $\rho^*$.
        \If{$\rho^*$ certifies $\pi$ is optimal} set $g(\pi)\leftarrow \boldsymbol{0}$.
        \Else 
            \State Select flows $\tilde{x}$ (e.g., that minimizes $\lvert\lvert \cdot\rvert\rvert_1$ node imbalances over $[\mathbf{\pi}]$).
            \State Compute $g(\pi)=c-A\tilde{x}$.
        \EndIf
    \EndIf
    \State \textbf{return} $g(\pi)$.
\end{algorithmic}
\end{algorithm}
In practice, $\pi^k\in \mathcal{T}_{\min}(\nset)\cup \mathcal{T}_{\max}(\nset)$ rarely, and when it happens, the subset of active flow variables whose flow is not uniquely determined by our rule is small in comparison to the overall number of flow variables. Thus, with judicious implementation, \cref{alg:improved_fw_gd} can inform optimal stopping criterion without significantly impacting overall efficiency.

Given the vast literature on convex optimization, network flows, optimal transport, etc., many other iterative approaches to solving $\fw$ can be utilized, and we make no attempt to detail them all here. In fact, our methods are clearly sub-optimal in the sense that they make no attempt to apply many of the well-known techniques shown to improve first-order methods (e.g. acceleration/momentum, regularization, or annealing). Additionally, our subgradient algorithm makes no use of the extremely structured nature of the underlying network (in contrast with the combinatorial approaches mentioned earlier). Since our primary objective was in motivating the geometric intuition behind any such iterative approach, we defer such improvements to future work.

\section{Discussion}\label{sec:conclusion}

In this paper we show that the tropical Fermat--Weber problem is dual to a min-cost flow problem, and that it can be formulated as an optimal transport problem for any positively weighted sample. We show that the set of optimal solutions is bi-tropically convex, and that it corresponds to a bounded cell of a particular polyhedral decomposition of $\TPT{n}$. This decomposition is generated by min- and max-plus tropical hyperplane arrangements whose apices correspond to points of the sample. This allows us to leverage existing algorithms and results on polytropes to efficiently obtain the set of all optimal tropical extreme points.

Tropical Fermat--Weber points are of particular interest in phylogenetics, where the underlying metric structure is well suited to comparing objects such as cophenetic vectors of coalescent trees. These vectors live in the \emph{space of ultrametrics} $\mathcal{U}$, a non-Euclidean subset of $\TPT{n}$ governed by a stronger version of the triangle inequality $\D(x,z)\leq \max(\D(x,y),\,\D(y,z))$. \cite{Aliatimis2024} showed that under their proposed model, a symmetric Fermat--Weber point in $\fw(\nset) \cap \mathcal{U}$ is an estimated species tree for a set of gene trees $\nset = \{\sample\} \subset \mathcal{U}$. As show in \citep{LSTY}, however, even if $\nset = \{\sample\} \subset \mathcal{U}$, $\fw(\nset) \not \subset \mathcal{U}$. Consequently, methods that compute or describe the full set of tropical Fermat--Weber points $\fw(\nset) \cap \mathcal{U}$ are especially relevant in this context. Thus, our work contributes an significant first step in this direction.

Throughout this work we assumed finite sample points $p_i \in \R^n$. Many of the tropical structures and techniques considered here extend naturally to the more general setting $p_i \in \Rmin^n$, where some coordinates may be non-finite. For instance, tropical hyperplanes require only two finite coefficients in their defining linear form. Extending the results and algorithms developed here to accommodate non-finite data therefore appears feasible and warrants further investigation.

Our other principal assumption concerned $m \ge n$ for the sample data; the case $m = n$ appears to be particularly revealing. Consider a discrete Wasserstein optimal transport problem specified by a cost matrix $C \in \Rmin^{m \times m}$ together with marginal vectors $a,b \in \R_{\ge 0}^m$ satisfying $\sum_i a_i = \sum_i b_i$. In the classical Wasserstein setting, the ground cost is induced by a metric and therefore satisfies the triangle inequality. In this case, $C$ coincides with its Kleene star and admits the factorization $C = -C^{\top} \odot C$.

Formulating the tropical Fermat--Weber problem in $\ot$ with $Q = C$, $\alpha = a$, and $\beta = b$ recovers precisely the original optimal transport linear program. Thus, every such Wasserstein optimal transport problem admits a (trivial) tropical Fermat--Weber representation. A more intriguing question is whether one can identify a matrix $V \in \R^{m \times n}$ with $n < m$ such that 
\[
-Q^{\top} \odot Q \approx C,
\]
with $n$ as small as possible. Such a representation would yield an approximation of the original transport problem via a min-cost flow formulation on a significantly sparser graph ($\mathcal{O}(nm)$ arcs rather than $\mathcal{O}(m^2)$). Investigating the existence, construction, and computational impact of such low-rank tropical factorizations may be worth exploring.

\section*{Acknowledgment}

RY, DB, and JS are partially supported from NSF DMS 1916037 and 2409819. Part of this research was performed while DB was visiting the Institute for Mathematical and Statistical Innovation (IMSI), which is supported by the National Science Foundation (Grant No. DMS-1929348).
KM is partially supported by JSPS KAKENHI Grant Numbers JP22K19816, JP22H02364 and JP25K15283. The views expressed in this article are those of the author(s) and do not reflect the official policy or position of the U.S. Naval Academy, Department of the Navy, the Department of War, or the U.S. Government.

\newpage

\appendix







\bibliography{refs}

\end{document}